\documentclass[a4paper,12pt]{article}
\usepackage{essentials}

\NewDocumentEnvironment{condition}{O{} m}
{
	\ifx&#1&%
	\def\criterionlabel{Condition\ (\textbf{#2})}%
	\else
	\def\criterionlabel{Condition\ \optionaldesc{(#2)}{#1}}%
	\fi
	\par\medskip
	\noindent\textbf{\criterionlabel.}\itshape
}
{
	\par\medskip
}

\title{A Criterion for Categories on which every\\Grothendieck Topology is Rigid}
\author{Jérémie Marquès}
\date{}

\bibliography{references.bib}

\begin{document}
\maketitle

\begin{abstract}
	Let $\cC$ be a small Cauchy-complete category. The subtoposes of $[\cC^\op,\cSet]$ are sometimes all of the form $[\cD^\op,\cSet]$ where $\cD$ is a full subcategory of $\cC$. This is the case for instance when $\cC$ is finite, an Artinian poset, or the simplex category. In order to unify these situations, we characterize the small categories $\cC$ such that for every $X ∈ \cC$, every subtopos of $[\cC_{/X}^\op,\cSet]$ is induced by a subcategory of $\cC_{/X}$. We provide two equivalent characterizations. The first one uses a two-player game, and the second one combines two ``local'' properties of $\cC$ involving respectively the poset reflections of its slices and its endomorphism monoids.
\end{abstract}

\begin{tcolorbox}[emphbox]
	\textbf{Disclaimer.} The open access published version \cite{MarCriterionCategoriesWhich2025} should be consulted instead. The current version incorrectly claims that every universally rigid category is Cauchy complete (the proof is wrong, I do not know whether this holds). This was pointed out by the reviewer and corrected in the published version. I take the opportunity to point out a similar result that appeared recently in \cite[Thm.~1.4]{DiLiLiaTorsionTheoreticInterpretation2025}, which I noticed only after publication.
\end{tcolorbox}

Let $\cC$ be a small Cauchy-complete category. If $\cD ⊆ \cC$ is a full subcategory, then $[\cD^\op,\cSet]$ is a subtopos of $[\cC^\op,\cSet]$. The Grothendieck topologies corresponding to these subtoposes are called \emph{rigid}. In this paper, we focus on the categories $\cC$ on which every Grothendieck topology is rigid. We call such a category \emph{universally rigid}. Several criteria ensuring universal rigidity can be found in the literature:
\begin{enumerate}[label=(C\arabic*)]
	✦ When the slices of $\cC$ are finite \cite[C2.2.21]{SketchesElephantToposJoh2002} \cite[Proposition~4.10]{CompleteLatticeEssentialKelLaw1989}.
	✦ When every arrow factors as a split epimorphism followed by a split monomorphism, and when every object has only a finite number of subobjects \cite[Cor.~3.3]{MenSuccessiveDimensionElegance2024}. \label{crit-2}
	✦ When $\cC$ is an Artinian poset \cite[Thm.~2.12]{GrothendieckTopologiesPosetLin2014}.
\end{enumerate}

This paper presents a criterion that generalizes all of the above. It falls short of being an exact characterization of universal rigidity, but it is modulo taking slices: We will characterize the small categories $\cC$ such that for every $X ∈ \cC$, every Grothendieck topology on $\cC_{/X}$ is rigid. We call such a category \emph{stably universally rigid}.

A mild necessary condition for $\cC$ to be universally rigid is that it is Cauchy complete. Otherwise, $[\cC^\op,\cSet] ≅ [\ovl{\cC}^\op,\cSet]$ where $\ovl{\cC}$ is the Cauchy-completion of $\cC$, and there is $X ∈ \ovl{\cC}$ not in $\cC$. Then the closure of $\set{X}$ under splitting of idempotents in $\ovl{\cC}$ yields a subtopos of $[\cC^\op,\cSet]$ not corresponding to a rigid topology on $\cC$. However, Cauchy-completeness is not sufficient. For instance, \cite[Thm.~2.12]{GrothendieckTopologiesPosetLin2014} shows that any non-Artinian poset admits a non-rigid Grothendieck topology.

\begin{rmq*}{}
	Neither \cite[Cor.~3.3]{MenSuccessiveDimensionElegance2024} nor the result presented here is more general. Indeed, \cite{MenSuccessiveDimensionElegance2024} is interested in the \emph{essential} subtoposes of a presheaf topos and \ref{crit-2} is only a special case of \cite[Cor.~3.3]{MenSuccessiveDimensionElegance2024}; the finiteness condition simply ensures that every subtopos is essential. There seems nonetheless to be a common theme, as the split epimorphisms also play an important role in the current paper.
\end{rmq*}

\paragraph{Acknowledgments} I thank Matías Menni for his comments, his advice to reorganize the paper, and for pointing out the connection with \cite{MenMonicSkeletaBoundaries2019,MenSuccessiveDimensionElegance2024}.

\paragraph{Funding} This research has been supported financially by the European Research Council (ERC) under the European Union's Horizon 2020 research and innovation program, grant agreement \#670624.

\paragraph{Notations} The composition of $f : X→Y$ and $g : Y→Z$ will be denoted $fg$.

\section{Three characterizations of stable universal rigidity}

In this section, we present three conditions \ref{SUR1}, \ref{SUR2} and \ref{SUR3} on a category $\cC$. The main theorem is that these conditions are all equivalent to stable universal rigidity.

We will consider a small category $\cC$ as sitting implicitly in $[\cC^\op,\cSet]$ via the Yoneda embedding. By the \emph{image} $\im(f)$ of a morphism $f$ in $\cC$, we mean the image of the corresponding morphism in $[\cC^\op,\cSet]$. With this convention, the morphisms in $\cC$ which are ``surjective'' are the split epimorphisms. The poset of subobjects of $X ∈ \cC$ of the form $\im(f)$ with $f : Y→X$ in $\cC$ can be identified with the \emph{poset reflection} of $\cC_{/X}$. The poset reflection of a small category is obtained by declaring that $A≤B$ when there is a morphism $A→B$, and by quotienting by the induced equivalence relation. A poset is \emph{Artinian} if there is no strictly decreasing sequence of elements.

\begin{condition}[SUR1]{SUR1}
\begin{enumerate}[label=(\roman*)]
	✦ The poset reflection of $\cC_{/X}$ is Artinian, for each $X ∈ \cC$. \label{SUR1-1}
	✦ For each $f : Y→X$, there is $f' : Y'→X$ with $\im(f) = \im(f')$ and such that for every $r : Z→Y'$, if $\im(f') = \im(rf')$ then $\im(r) = Y'$. \label{SUR1-2}
\end{enumerate}
\end{condition}

The \emph{game of split epi} on a small category $\cC$ is defined as follows. There are two players: \emph{Cleaner} and \emph{Reducer}. The state of the game is determined by a morphism $f$ in $\cC$ whose codomain never changes. A play by Cleaner in position $f$ consists in choosing a factorization $f = rf'$ where $r$ is a split epimorphism; the state of the game changes from $f$ to $f'$. Dually, a play by Reducer in position $f$ consists in writing $f' = rf$ where $r$ is \emph{not} a split epimorphism; the state of the game changes from $f$ to $f'$. Cleaner wins if Reducer cannot make a valid move, and Reducer wins if the game continues forever.

% We could allow more generally Cleaner to factorize $f$ as $rf'$ with $r$ arbitrary. However, if $r$ is not a split epimorphism, this move is useless since Reducer can just cancel it by precomposing by $r$.

Alternatively, a move by Cleaner when the game is in position $f : Y→X$ is determined by an idempotent $e : Y→Y$ such that $ef = f$. We split $e = rs$ and the position of the game becomes $sf$. This corresponds to the factorization $f = r(sf)$. Note that Cleaner cannot precompose $f$ by an \emph{arbitrary} split monomorphism, because of the condition $ef=f$.

Cleaner never modifies the image of $f$ and Reducer always makes it smaller, although not necessarily strictly smaller. We think of a game position $f : Y→X$ as a ``presentation'' of $\im(f) ⊆ X$. The goal of Reducer is to continuously make $\im(f)$ smaller, or at least \emph{pretend to}. The goal of Cleaner is to ``clean up'' that presentation, so as to prevent Reducer from lying indefinitely.

\begin{condition}[SUR2]{SUR2}
	Cleaner has a winning strategy in the game of split epi.
\end{condition}

We say that a monoid $M$ has \emph{enough idempotents on the left} if for every $f ∈ M$, either $f$ has a left inverse or there is a non identity idempotent $e$ and some $n≥1$ satisfying $ef^n = f^n$. Every finite monoid has enough idempotents on the left, by taking $e = f^ω$ the only idempotent power of $f$.

\begin{condition}[SUR3]{SUR3}
\begin{enumerate}[label=(\roman*)]
	✦ $\cC$ is Cauchy-complete.\label{SUR3-1}
	✦ The poset reflection of $\cC_{/X}$ is Artinian, for each $X ∈ \cC$.\label{SUR3-2}
	✦ $\End(X)$ has enough idempotents on the left, for each $X ∈ \cC$.\label{SUR3-3}
\end{enumerate}
\end{condition}

\begin{thm}{}{main}
	The conditions \ref{SUR1}, \ref{SUR2} and \ref{SUR3} are all equivalent to stable universal rigidity.
\end{thm}

As an example, we show Corollary~\ref{cor:degree-criterion} that applies to the simplex category.

\begin{cor}{}{degree-criterion}
	Let $\cC$ be a small category equipped with an ordinal-valued degree function $d : \abs{\cC}→\Ord$ such that every morphism factors as a split epimorphism followed by either an isomorphism or a morphism $f : A→B$ with $d(A) < d(B)$. Then $\cC$ is universally rigid.
\end{cor}

\begin{proof*}{}
	We show \ref{SUR2}. The winning strategy of Cleaner is to react to the moves of Reducer as follows. Suppose that the previous move of Reducer was to precompose $f : Y→X$ by $h : Z→Y$. Then Cleaner factorizes $h = rg$ with $r$ the split epimorphism given by the hypothesis on the degree function.
	\[\begin{tikzcd}
		Z \ar[r,"h"] \ar[d,->>,"r"] & Y \ar[r,"f"] & X\\
		Y' \ar[ru,"g"']
	\end{tikzcd}\]
	The morphism $g$ cannot be invertible as it would imply that $h$ is a split epimorphism. Hence $d(Y') < d(Y)$. As the degree strictly decreases, the game cannot continue forever.
\end{proof*}

\section{The proof}
\label{sec:proof}

We split the proof of Theorem~\ref{thm:main} in a sequence of lemmas.

We will use the following definition from \cite[C2.2.18]{SketchesElephantToposJoh2002}. Given a Grothendieck topology $J$ on $\cC$, an object $X ∈ \cC$ is called \emph{$J$-irreducible} if the only covering sieve is the maximal sieve. Then $J$ is rigid if and only if every object admits a $J$-covering by $J$-irreducible objects. When $J$ is rigid, the subtopos $\Sh(\cC,J)$ is equivalent to $[\cD^\op,\cSet]$ where $\cD$ is the full category of $\cC$ spanned by the $J$-irreducible objects.

\begin{lem}{}{sufficient}
	If $\cC$ satisfies \ref{SUR1}, then $\cC$ is universally rigid.
\end{lem}

\begin{proof*}{}
	Suppose that $\cC$ satisfies \ref{SUR1}. Let $J$ be a Grothendieck topology on $\cC$ and let $X ∈ \cC$. We show that for all $f : Y→X$, the image $\im(f)$ of $f$ is $J$-covered by $J$-irreducible objects of $\cC$. We proceed by Artinian induction on $\im(f)$, thanks to \ref{SUR1}\ref{SUR1-1}. Taking $f$ the identity of $X$, this will show that $X$ is $J$-covered by $J$-irreducible objects.
		
	Let $f : Y→X$ and suppose that the target property holds for all $g : W→X$ with $\im(g) < \im(f)$. Thanks to \ref{SUR1}\ref{SUR1-2}, we can suppose that for all $r : Z→Y$, if $\im(rf) = \im(f)$, then $\im(r) = Y$. If $Y$ is $J$-irreducible, then $\im(f)$ is covered by $Y$ and we are done. Otherwise, there is some $J$-covering $S$ of $Y$ which does not generate the maximal sieve, hence which contains no split epimorphism. We thus have $\im(rf) < \im(f)$ for every $r ∈ S$, and by induction there is a $J$-covering of $\im(rf)$ consisting of $J$-irreducible objects. This transfers to $\im(f)$ by composition of $J$-coverings.
\end{proof*}

\begin{rmq}{}{Cauchy-complete}
	Condition \ref{SUR1} implies that $\cC$ is Cauchy complete, since it is necessary in order for $\cC$ to be universally rigid. Cauchy-completeness is actually a special case of \ref{SUR1}\ref{SUR1-2}, when $f$ is idempotent.
\end{rmq}

\begin{lem}{}{game}
	The conditions \ref{SUR1}, \ref{SUR2} and \ref{SUR3} are equivalent.
\end{lem}

\begin{proof*}{}
	\proofstep{\ref{SUR1} $⇒$ \ref{SUR2}}
	Suppose that \ref{SUR1} holds. The winning strategy of Cleaner in position $f : Y→X$ is to answer with the morphism $f' : Y'→X$ given by \ref{SUR1}\ref{SUR1-2}. Since $\im(f) = \im(f')$, there is a morphism $r : Y'→Y$ such that $f = rf'$, and $r$ is a split epimorphism because we must have $\im(r) = Y'$. In the next move, Reducer is forced to produce a morphism $f'' : Y''→X$ with $\im(f'') < \im(f')$. By \ref{SUR1}\ref{SUR1-1}, the process must terminate.
	
	\proofstep{\ref{SUR2} $⇒$ \ref{SUR1}} Suppose that the poset reflection of $\cC_{/X}$ is not Artinian for some $X ∈ \cC$. We produce a winning strategy for Reducer for the game starting in position $\id_X : X→X$. Let $(f_i : Y_i → X)_{i∈ℕ}$ be such that $\im(f_{i+1}) < \im(f_i)$ for all $i ∈ ℕ$. Then Reducer can win by playing $f_i$ as its $i$th move. Cleaner answers with $f'_i : Y'_i → X$ such that $\im(f'_i) = \im(f_i)$. Since $\im(f_{i+1}) < \im(f_i) = \im(f'_i)$, there is an arrow $r : Y_{i+1} → Y'_i$ with $rf'_i = f_{i+1}$. The inclusion being strict, $r$ is not a split epimorphism. Hence Reducer can play $f_{i+1}$ as its next move, and so on.
	
	Now, suppose that \ref{SUR1}\ref{SUR1-2} does not hold. This means that there is some $f : Y→X$ such that for every $f' : Y'→X$ satisfying $\im(f) = \im(f')$, there is an $r : Z→Y'$ such that $\im(r) < Y'$ but $\im(rf') = \im(f')$. We will again produce a winning strategy for Reducer starting from the identity $\id_X : X→X$. The first move of Reducer is used to put the game in position $f : Y→X$. This is possible since $\im(f) < X$, otherwise we could take $f' = \id_X$. Afterward, the game will always be in a position $f' : Y'→X$ such that $\im(f') = \im(f)$. Cleaner cannot change this. Reducer can then use the $r : Z→Y'$ given by the hypothesis and put the game in state $rf'$. This is a valid move since $r$ is not a split epimorphism.
	
	\proofstep{\ref{SUR3} $⇒$ \ref{SUR1}} Suppose that $\cC$ satisfies \ref{SUR3}. We must show \ref{SUR1}\ref{SUR1-2}. Let $f : Y→X$. If there is a non-identity idempotent $e_1 : Y→Y$ such that $e_1 f = f$, we split $e_1$ as $r_1s_1$ and we define $f_1 = s_1f$. Then $\im(f) = \im(f_1)$, because $\im(f) = \im(r_1s_1f) ⊆ \im(s_1f) = \im(f_1)$. Repeat the process as long as there is a non-identity idempotent $e_{i+1}$ with $e_{i+1}f_i = f_i$. The process must terminate, otherwise $\im(s_1), \im(s_2s_1), …$ is an infinite strictly decreasing sequence of subobjects of $X$. Let $f' : Y'→X$ be the morphism obtained at the last step, so that the only idempotent $e : Y'→Y'$ verifying $ef' = f'$ is the identity. Let $r : Z→Y'$ such that $\im(rf') = \im(f')$. This means that there is some $s : Y'→Z$ such that $srf' = f'$. We wish to show that $\im(r) = Y'$, i.e., that $r$ has a left inverse. Since $\End(Y')$ has enough idempotents on the left, there are two possibilities. Suppose there is a non identity idempotent $e : Y'→Y'$ and $n≥1$ such that $e(sr)^n = (sr)^n$. That would imply that $ef' = e(sr)^nf' = (sr)^nf' = f'$, but this is impossible by the hypothesis on $f'$. Hence $sr$ has a left inverse, and $r$ too.
	% We actually use that the idempotents with the Mitsch partial order embed in the poset reflection of $\cC_{/X}$.
	
	\proofstep{\ref{SUR1} $⇒$ \ref{SUR3}} Suppose that $\cC$ satisfies \ref{SUR1}. Then $\cC$ is Cauchy-complete by Remark~\ref{rmq:Cauchy-complete}. We show that if \ref{SUR3}\ref{SUR3-3} does not hold, then Reducer has a winning strategy, contradicting \ref{SUR2}. Suppose that $f : X→X$ is an endomorphism in $\cC$ with no left inverse and such that there is no non-identity idempotent $e : X→X$ such that $ef^n = f^n$ for $n≥1$. The game starts in position $f : X→X$. The winning strategy of Reducer is to always play $f$. The successive positions of the game are $f$, $f^2$, $f^3$, etc. At each step, Cleaner can only play the identity since it is the only idempotent $e$ such that $ef^n = f^n$, for $n≥1$.
\end{proof*}

\begin{lem}{}{slice}
	The equivalent conditions \ref{SUR1}, \ref{SUR2} and \ref{SUR3} are stable under taking slices.
\end{lem}

\begin{proof*}{}
	We use \ref{SUR2}. Note that the game is played in a slice $\cC_{/X}$, except that the valid moves of Cleaner and Reducer depend on which morphisms are split epi in $\cC$, and not in $\cC_{/X}$. Nonetheless, we show that they coincide with the split epi in $\cC_{/X}$. Let $f : Y→X$ and $f' : Y'→X$ be two objects in $\cC_{/X}$, and let $g : Y'→Y$ be a morphism in over $X$ which is a split epi in $\cC$. This means that there is $g' : Y→Y'$ such that $g'g = \id_Y$. But $g'f' = g'gf = f$, so that $g'$ is a morphism in $\cC_{/X}$ and $g$ is a split epi in $\cC_{/X}$.
\end{proof*}

The last step of the proof of Theorem~\ref{thm:main} involves building a non-rigid topology in case $\cC$ does not satisfy \ref{SUR3}. Just like in \cite{GrothendieckTopologiesPosetLin2014}, the main tool for this will be the \emph{double negation topology}. This topology will be denoted by $J_{¬¬}$, leaving the category implicit. A sieve on $X$ is $J_{¬¬}$-dense when it has a non-empty intersection with every non-empty sieve on $X$.

\begin{lem}{}{negneg-irred}
	An object $X$ of a small category is $J_{¬¬}$-irreducible if and only if every morphism $Y→X$ is a split epimorphism.
\end{lem}

\begin{proof*}{}
	If every morphism $Y→X$ is split, then any non-empty $J_{¬¬}$-covering sieve on $X$ contains the identity, hence $X$ is $J_{¬¬}$-rigid.
	
	Reciprocally, suppose that there is $f : Y→X$ which is not a split epimorphism. Let $C$ be the sieve on $X$ consisting of all the morphisms which are not split epimorphisms. For any $g : Z→X$, either $g$ is in $C$, or it is a split epimorphism. In the latter case, there is a section $s : X→Z$ and $f = (fs)g$. This shows that $C$ is covering with respect to the double negation topology, despite not being maximal, hence $X$ is not $J_{¬¬}$-irreducible.
\end{proof*}

\begin{lem}{}{converse}
	If $\cC$ is stably universally rigid, then it satisfies \ref{SUR3}.
\end{lem}

\begin{proof*}{}
	We will show that if $\cC$ does not satisfy \ref{SUR3}, then it is not stably universally rigid. As already seen, if $\cC$ is not Cauchy-complete, then it is not universally rigid.
	
	Suppose that the poset reflection of $\cC$ is not Artinian. This means that there is a chain $X_0 ← X_1 ← ⋯$ with no morphism $X_i → X_j$ when $j > i$. Let $\cD$ be the Cauchy completion of the full subcategory $\set{X_0,X_1,…} ⊆ \cC$. Let $J_{¬¬}$ be the double negation topology on $\cD$. By Lemma~\ref{lem:negneg-irred}, there is no $J_{¬¬}$-irreducible object in $\cD$. Since it is not empty, the subtopos $\Sh(\cD,J_{¬¬}) ⊆ [\cD^\op,\cSet]$ does not come from a subcategory of $\cD$, hence the subtopos $\Sh(\cD,J_{¬¬}) ⊆ [\cC^\op,\cSet]$ does not come from a subcategory of $\cC$.
	
	If $\cC_{/X}$ is not Artinian for some $X ∈ \cC$, the previous argument applies and $\cC$ is not stably universally rigid.
	
	Suppose that $\End(X)$ does not have enough idempotents on the left for some $X ∈ \cC$. Let $f : X→X$ such that $f$ is not left invertible and such that for all $n$, there is no non-identity idempotent $e : X→X$ such that $ef^n = f^n$. Let $\cD$ be the full subcategory $\set{f^0,f^1,f^2,…} ⊆ \cC_{∕X}$. By the assumption, $\cD$ is Cauchy-complete because there is no non-trivial idempotent. Moreover, for every $n ∈ ℕ$, the morphism $f : f^{n+1} → f^n$ is not a split epimorphism, hence by Lemma~\ref{lem:negneg-irred} there is no irreducible object relatively to the double negation topology. Once again, this shows that $\Sh(\cD,J_{¬¬}) ⊆ [\cC_{∕X}^\op,\cSet]$ is a subtopos which does not correspond to a subcategory of $\cC_{/X}$.
\end{proof*}

\section{Concluding remarks}

One could wonder whether the conditions \ref{SUR3}\ref{SUR3-1} and \ref{SUR3}\ref{SUR3-2} are enough to deduce that $\cC$ is universally rigid. Let us explain why this is a natural question. The functor $X ↦ \cPoset(\cC_{∕X})$ sending each object $X$ to the poset reflection of $\cC_{∕X}$ is a \emph{polyadic poset} in the sense of \cite[§~2.2, p.~69]{MarCategoricalLogicPerspective2023}. Viewing each poset $P$ as an Alexandroff locale, the functor  $X ↦ \cPoset(\cC_{∕X})$ is also a polyadic locale, more commonly called a geometric hyperdoctrine \cite[Ch.~III]{WriDoctrinalGroupoidalRepresentations2024}. The classifying topos of this polyadic locale is $[\cC^\op,\cSet]$, so the subtoposes of $[\cC^\op,\cSet]$ are in bijection with its polyadic sublocales. If $\cPoset(\cC_{∕X})$ is Artinian for each $X$, we know by \cite[Thm.~2.12]{GrothendieckTopologiesPosetLin2014} that each polyadic sublocale is actually a polyadic subposet. This seems to be an important step towards universally rigidity, but the following counter-example shows that this is not enough. It would nonetheless be interesting to understand how to combine this remark with \ref{SUR3}\ref{SUR3-3} to deduce universal rigidity.

\begin{exmp}{}{}
	Let $S$ be a non-empty semigroup with no idempotent and only one left ideal, meaning that $Sa = S$ for all $a ∈ S$. For instance, $S$ can be the opposite of the semigroup of all the injections $f : ℕ→ℕ$ such that $\card{ℕ⧵f[ℕ]} = \card{ℕ}$. Let $M = S∪\set{1}$ be the monoid obtained by adding an identity to $S$. Then $M$ seen as a one-object category satisfies \ref{SUR3}\ref{SUR3-1} and \ref{SUR3}\ref{SUR3-2}. But by Lemma~\ref{lem:negneg-irred}, there is no irreducible object for the double negation topology on $M$, hence this topology is not rigid.
\end{exmp}

We defined stable universal rigidity by requiring that $\cC_{/X}$ be universally rigid for every $X ∈ \cC$. In fact, stable universal rigidity satisfies a stronger stability property: $X$ can be replaced by any presheaf. When $F : \cC^\op→\cSet$ is a presheaf, we denote by $\cC_{/F}$ the category of elements of $F$. We claim that stable universal rigidity transfers from $\cC$ to $\cC_{/F}$. To see that, we use that $(\cC_{/F})_{/(x,X)} ≅ \cC_{/X}$ for any $x ∈ F(X)$, plus the fact that \ref{SUR2} only depends on the slices of $\cC$ (see Lemma~\ref{lem:slice}). In light of the equivalence $[\cC_{∕F}^\op,\cSet] ≅ [\cC^\op,\cSet]_{∕F}$, if we identify a Cauchy-complete category with the corresponding presheaf topos, this means that stable universal rigidity is not only stable under taking slices over representables, but over any presheaf. Thus, stable universal rigidity ensures that every subtopos of a slice of $[\cC^\op,\cSet]$ is of a certain simple form, in particular a presheaf topos. In terms of the theory of internal locales of \cite{JoyTieExtensionGaloisTheory1984}, the subtoposes of the slices of a topos correspond to the sublocales of the discrete internal locales, and the corresponding geometric morphisms are stable under composition. These remarks might hint at a more conceptual understanding of stable universal stability, even though the initial motivation is mostly technical.

\AtNextBibliography{\small}
{\printbibliography[
heading=bibintoc,
title={References}
]}

\end{document}